\documentclass[12pt]{article}
\usepackage{amsmath}
\usepackage{amssymb}
\usepackage{amsthm}
\usepackage{mathtools}

\usepackage{enumitem}

\newcommand{\dt}{\partial_t}
\newcommand{\dx}{\partial_x} 
\newcommand{\idx}{\,dx}

\newcommand{\norm}[2]{\Vert #1 \Vert_{#2}}

\newcommand{\initial}{\text{initial data}}
\newcommand{\mrm}[1]{\mathrm{#1}}

\newtheorem{theorem}{Theorem}

\theoremstyle{remark}
\newtheorem{remark}{Remark}

\numberwithin{equation}{section}

\title{Viscous pressureless flows with free boundary in one space dimension: The constant viscosity case}
\author{Xin Liu\footnote{Texas A\&M University, College Station, TX. xliu23@tamu.edu}}

\begin{document}

\maketitle

\begin{abstract}
    We establish the global well-posedness of the free boundary problem of the viscous pressureless and almost pressureless heat conductive flows in one space dimension. In both cases, arbitrarily large but smooth initial data is considered, and the evolving fluid domains remain bounded for all time. In the viscous pressureless case, we are able to identify the terminal flow domain in terms of the initial data. In the viscous almost pressureless case, we construct the flow as a perturbation of the viscous pressureless flow, and establish the first result for the Navier-Stokes-Fourier system in the current setting.


\end{abstract}


\section{Introduction}

\subsection{Backgrounds and the scope of this paper}

The pressureless flow model \eqref{sys:pressureless-flows} is used to describe dynamics of sticky particle \cite{bouchutZeroPressureGas1994,brenierStickyParticlesScalar1998}, where two particles stick to each other with the same terminal velocity after collision. It is an extension (or a simplification, respectively) of Burgers' equation (or the compressible Euler equations, respectively). In the spirit of \cite{dipernaOrdinaryDifferentialEquations1989}, the pressureless Euler equations are strongly tied to dynamical systems with large dimension, describing the motion of large number of particles. Indeed, it can be also used in astrophysics to model the mass distribution of large scale structures of the galaxy \cite{vergassolaBurgersEquationDevils1994,zeldovichGravitationalInstabilityApproximate1970}. It is also related to the modeling of traffic flow \cite{zatorskaRecentAdvancesAnalysis2025}. 

In comparison to the formation of singularity in Burgers's equation and the compressible Euler equations, new singularity is observed in the pressureless Euler equations. By investigating the Riemann problem in one-dimension, the authors in \cite{chenFormation$delta$ShocksVacuum2003,chenConcentrationCavitationVanishing2004,shenFormationDeltaShocks2010,shengConcentrationMassPressureless2020} have rigorously shown the formation of $ \delta $-shock and vacuum states in the pressureless limit from the compressible Euler equations. To be more precise, for the pressureless Euler equations, a weighted $ \delta $ measure will form from a compressing initial data, and vacuum, instead of rarefaction wave, will form from a depressing initial data. In the case of multi-dimensional flows, the problem of non-uniqueness and non-existence arise, as shown in \cite{bressanNonexistenceNonuniquenessMultidimensional2014,bressanGenericSingularities2D2025}.

Meanwhile, the measure value weak solutions to the pressureless Euler equations have been intensively studied. For instance, a generalized variational principle in the spirit of the Lax-Oleinik variational principle for scalar conservation laws, was established to construct global weak solutions to the pressureless Euler equations \cite{eGeneralizedVariationalPrinciples1996}. Similar result can be found also in \cite{wangCauchyProblemTransportation1997}. The uniqueness of the variational weak solutions is shown in \cite{huangWellPosednessPressureless2001}. The author in \cite{boudinSolutionBoundedExpansion2000} constructs unique global smooth solutions to the viscous pressureless flow, and establishes the unique weak solution in the sense of duality, introduced in \cite{bouchutOnedimensionalTransportEquations1998,bouchutDualitySolutionsPressureless1999}, as the vanishing viscosity limit. We refer readers to \cite{brenierStickyParticlesScalar1998,huangWeakSolutionPressureless2005,moutsingaConvexHullsSticky2008,natileWassersteinApproachOneDimensional2009,ducometBoundaryValueProblem2009,breschSingularLimitNavier2014,nguyenOneDimensionalPressurelessGas2015,cavallettiSimpleProofGlobal2015,figalliRigorousDerivationKinetic2018,danilovInteractionDShockWaves2019,carrilloPressurelessEulerNonlocal2020,stefanoStickyParticleSolutions2020,neumannInitialboundaryValueProblem2022,klyushnevModelTwoDimensionalPressureless2023,heExactSphericallySymmetric2024,carrilloEquivalenceEntropySolutions2024,dongAnalyticalSolutionsPressureless2024,sen1DPressurelessGas2025,learGeometricStructureMass2022,carrilloQuantifyingHydrodynamicLimit2021,liuStabilityDeltaWave2022,liangVanishingPressureLimit2017} and the references therein for other construction of solutions and recent works on the pressureless flows. 

The goal of this paper is to explore the effect of viscosity in the pressureless flow. To be more precise, by ignoring the pressure, we aim at isolating the effect of the viscosity, describing the trajectory of smooth solutions, and establish the connection of the viscous pressureless flow to the Navier-Stokes-Fourier flow. In the one-dimensional free boundary setting, where the fluid boundary is moving together with the flow, we show first, for arbitrary smooth but large initial data, the viscous pressureless flow \eqref{sys:pressureless-flows} has a unique global smooth solution, and the domain of the flow exponentially converges to a compact set as time goes to infinity. Moreover, we are able to characterize the terminal flow domain in terms of the initial data only. We emphasize that, this asymptotic behavior is the result of the constant viscosity, in comparison to the inviscid flow and the density-dependent viscosity flow \cite{heExactSphericallySymmetric2024}.

In addition, in the hypersonic regime, i.e., almost pressureless case, we show that for arbitrary smooth but large initial data, there exists a unique global solution to the free boundary problem of the compressible Navier-Stokes-Fourier system \eqref{sys:NSF-flow} with dirichlet boundary condition for the temperature. The resulting flow should be seen as a perturbation of the pressureless flow and as a result, the flow domain remains a compact set as time goes to infinity. We emphasize that this result is not possible in the setting of zero normal heat flux, where the latter is establish in \cite{chenGLOBALSOLUTIONSNAVIERSTOKES2002}. This is due to the difference in asymptotic behavior for solutions to the heat equation with dirichlet and zero normal heat flux boundary conditions. 

Our results are the first to tackle the free boundary problem of (almost) pressureless flow. The purpose of this study is to distinguish the physical effect of the viscosity to the asymptotic dynamics of (pressureless) fluid flow. 

We would like to make a few remarks before finishing the introduction. The lack of dispersive background flows (cf. \cite{LY2019,hadvzicExpandingLargeGlobal2018,hadvzicClassGlobalSolutions2019}) means that we don't have mechanism to prevent splash and splat singularities (e.g. \cite{Coutand2014}) in the case of multi-dimensional flow. Therefore our results cannot be trivially extended to higher dimension without any restrictions on geometry, initial data, etc.. The effect of degenerate viscosity, i.e., density and/or temperature depending viscosity, is dramatically different from the constant viscosity as shown in the appendix and \cite{heExactSphericallySymmetric2024}, and therefore deserves its separate study. We prefer to not complicating our current paper and leave these issues to future works.

\subsection{Main results}

Let $ \rho \geq 0 $ be the flow density, and $ u $ be the flow velocity. The viscous pressureless flow is governed by
\begin{subequations}\label{sys:pressureless-flows}
    \begin{gather}
        \label{eq:continuity} \partial_t \rho + \partial_y (\rho u) = 0, \\
        \label{eq:momentum} \partial_t(\rho u) + \partial_y (\rho u^2) = \partial_y (\mu \partial_y u).
    \end{gather}
    Here $ \mu \geq 0 $ is the viscosity coefficient, which is taken to be constant in this work. 

    We consider the flow in a simply connected, evolving domain, where the domain boundary moves together with the flow without surface stress. That is
    \begin{align}
        \text{domain of the flow:} & \quad  \Omega(t) = (a(t), b(t)), \label{def:domain} \\
        \text{evolution of the boundary:} & \quad  \dot a(t) = u(a(t), t), \quad \dot b(t) = u(b(t),t), \label{bc:moving-boundary} \\
        \text{stress free boundary condition:} & \quad \mu \partial_y u\vert_{y = a(t), b(t)} = 0. \label{bc:stress-free}
    \end{align}

\end{subequations}
Then we have the following theorem for the viscous pressureless flow:
\begin{theorem}[Pressureless flow]
    \label{thm:pressureless} Denote by $ (\rho,u)\vert_{t=0}=(\rho_0, u_0) $ the initial data and $ \rho_0\vert_{\Omega(0)} >0 $. Without loss of generality, assume that the total momentum of the initial data is zero, i.e., 
    \begin{equation}
        \label{eq:initial-momentum} \int_{\Omega(0)} \rho_0 u_0 \ dy = 0.
    \end{equation}
    Suppose the initial data is smooth and satisfying the following:
    \begin{equation}
        \label{initial_ene:pressureless}
        \begin{gathered}
            \int_{\Omega(0)} \rho_0 \, dy + \int_{\Omega(0)} \rho_0 \vert u_0 \vert^2 \ dy + \int_{\Omega(0)} \rho_0 \vert u_{1} \vert^2 \ dy \leq \mathfrak E_0 < \infty.
        \end{gathered}
    \end{equation}
    for some $ \mathcal E_0 \in (0,\infty) $. Here $ u_{1} := u_t\vert_{t=0} $ is given by the equation \eqref{eq:momentum} at $ t = 0 $. Then there exists a global strong solution $ (\rho,u) $ to system \eqref{sys:pressureless-flows} with 
    \begin{equation}
        \label{asmpt:pressureless}
        \vert \Omega(t) \vert = \vert a(t) - b(t) \vert \xrightarrow[]{t \rightarrow \infty} \int_{a(0)}^{b(0)} \exp\bigl( - \frac{1}{\mu} \int_{a(0)}^{y} \rho_0 u_0 \ dy' \bigr) \ dy,
    \end{equation}
    and
    \begin{equation}
        \label{ene:total-pressureless}
        \begin{gathered}
        \sup_{0\leq t < \infty} \biggl\lbrace \int_{\Omega(t)} \rho \vert u\vert^2 \,dy + \int_{\Omega(t)} \rho \vert u_t \vert^2 \ dy \biggr\rbrace \\
        + \int_0^\infty \biggl( \int_{\Omega(t)}   \vert u_y \vert^2  \,dy + \int_{\Omega(t)} \vert u_{yt} \vert^2 \ dy \biggr) dt < C \mathfrak E_0,
    \end{gathered}
    \end{equation}
    for some constant $ C \in (0,\infty) $. 

\end{theorem}
\begin{remark}
    \begin{enumerate}[label = (\roman*)]
        \item We emphasize that there is no smallness assumption of $ \mathcal E_0 $.
        \item The assumption \eqref{eq:initial-momentum} can be removed easily by considering the Galilean transformation. 
        \item The detailed energy decay estimate is deferred to the Lagrangian description, where it is easier to describe. 
    \end{enumerate}
\end{remark}

\bigskip 

On the other hand, take into account the temperature $ \theta \geq 0 $, and let $ M > 0 $ be a large parameter, to be determined later. We consider the following Navier-Stokes-Fourier system for heat conductive flow:
\begin{subequations}\label{sys:NSF-flow}
    \begin{gather}
        \label{eq:NSF-continuity} \partial_t \rho + \partial_y (\rho u) = 0, \\
        \label{eq:NSF-momentum} \partial_t(\rho u) + \partial_y (\rho u^2) + \dfrac{1}{M^2} \partial_y(\rho \theta) = \partial_y (\mu \partial_y u), \\
        \label{eq:NSF-thermodynamics} \partial_t(\rho \theta) + \partial_y(\rho \theta u) + \dfrac{1}{M^2} \rho \theta \partial_y u = \partial_y (\kappa \partial_y \theta ) + \mu \vert \partial_y u \vert^2, 
    \end{gather}
    with, in addition to the boundary conditions \eqref{def:domain} -- \eqref{bc:stress-free}, 
    \begin{equation}
        \label{bc:thermo} \theta\vert_{y = a(t), b(t)} = 0.
    \end{equation}
\end{subequations}
We have the following theorem for this almost pressureless system: 
\begin{theorem}[Almost pressureless flow]
    \label{thm:almost-pressureless} 
    In addition to the assumption in theorem \ref{thm:pressureless}, assume that
    \begin{equation}
        \int_{\Omega(0)} \rho_0 \theta_0\,dy + \int_{\Omega(0)} \rho_0 \vert \theta_0\vert^2 \,dy + \int_{\Omega(0)} \rho_0 \vert \theta_1 \vert^2 \,dy \leq \mathfrak E_1 < \infty,
    \end{equation}
    where $ \theta_0 := \theta\vert_{t=0} $ and $ \theta_1 := \theta_t\vert_{t=0} $, with the latter defined by the equation \eqref{eq:NSF-thermodynamics}. Then for $ M $ sufficiently large, depending on the initial data, there exists a unique global strong solution $ (\rho,u, \theta) $ to system \eqref{sys:NSF-flow} with $ \vert \Omega(t) \vert < C < \infty $ for all $ t \in [0,\infty) $, and the following estimates hold:
    \begin{equation}
        \label{ene:total-almost-pressureless}
        \begin{gathered}
            \sup_{0\leq t < \infty} \biggl\lbrace \int_{\Omega(t)} \rho \vert u\vert^2 \,dy + \int_{\Omega(t)} \rho \vert u_t \vert^2 \ dy + \int_{\Omega(t)} \rho \theta \,dy \\ + \int_{\Omega(t)} \rho \vert \theta \vert^2 \,dy + \int_{\Omega(t)} \rho \vert \theta_t \vert^2 \,dy \biggr\rbrace \\
            +  \int_0^\infty \biggl( \int_{\Omega(t)}   \vert u_y \vert^2  \,dy + \int_{\Omega(t)} \vert u_{yt} \vert^2 \ dy + \int_{\Omega(t)} \vert \theta_y \vert^2 \,dy + \int_{\Omega(t)} \vert \theta_{yt} \vert^2 \,dy \biggr) dt \\
             < C (\mathfrak E_0 + \mathfrak E_1 + \mathfrak E_0^3 + \mathfrak E_1^3),
        \end{gathered}
    \end{equation}
    for some constant $ C \in (0,\infty) $. 
\end{theorem}

\begin{remark}
    \begin{enumerate}[label = (\roman*)]
        \item System \eqref{sys:NSF-flow} can be obtained by rescaling the Navier-Stokes-Fourier system, and it preserves the mass, momentum, and energy conservation laws. 
        \item The result in theorem \ref{thm:almost-pressureless} should be seen as a perturbation of the result of theorem \ref{thm:pressureless}.
        \item By scaling the temperature $ \widetilde{\theta}:= \theta/M^2 $, the result in theorem \ref{thm:almost-pressureless} can be seen as a solution with arbitrary velocity but small rescaled temperature $ \widetilde{\theta} $. 
        \item Again, we refer to the Lagrangian description for the detailed energy decay estimates.
    \end{enumerate}
\end{remark}

For both theorems \ref{thm:pressureless} and \ref{thm:almost-pressureless}, we will only show the {\it a prior} estimates for the solution. To existence and uniqueness of solutions follow from the standard theory of partial differential equations \cite{Evens2010}, and are left out in order to simply the presentation.

The rest of the paper will be organized as follows. In section \ref{sec:pressureless-flow}, we focus on the viscous pressureless flow \eqref{sys:pressureless-flows} and aim at proving theorem \ref{thm:pressureless}. More precisely, in section \ref{subsec:la-fmlt}, we reformulate the system of equations in the Lagrangian coordinates. The most important ingredient, the point-wise estimate, will be given in section \ref{subsec:pointwise-pressureless}. Finally, in section \ref{subsec:regularity-pressureless}, we obtain the global regularity for the viscous pressureless flow and finish the proof of theorem \ref{thm:pressureless}. The proof of theorem \ref{thm:almost-pressureless} is given in section \ref{sec:almost-pressureless}. Following similar structure as before, we show the important point-wise estimate in section \ref{subsec:pointwise-nsf}. Then the {\it a prior} estimates are given in section \ref{subsec:estimates-nsf}. Finally, we conclude the proof of theorem \ref{thm:almost-pressureless} in section \ref{subsec:collecting-nsf}.

\section{Viscous pressureless flow with constant viscosity}
\label{sec:pressureless-flow}

\subsection{Lagrangian formulation}\label{subsec:la-fmlt}

We write system \eqref{sys:pressureless-flows} in the Lagrangian coordinates. Without loss of generality, assume that $ \Omega(0) = (-1,1) $. Let $ y = \eta(x,t) $ satisfy 
\begin{equation}\label{def:lagrangian}
    \partial_t \eta(x,t) = u(y,t) = u(\eta(x,t),t), \qquad \eta(x,0) = x \in (-1,1). 
\end{equation}
Meanwhile, write the density and the velocity in the Lagrangian coordinates $ (x,t) $ by
\begin{equation}\label{def:lg-variables}
    f(x,t) := \rho(\eta(x,t),t), \qquad v (x,t) := u(\eta(x,t),t).
\end{equation}
Then it is easy to verify from \eqref{eq:continuity} and \eqref{def:lagrangian} that $$ f(x,t) = \dfrac{\rho_0(x)}{\eta_x(x,t)}. $$

Therefore, system \eqref{sys:pressureless-flows} can be written in the lagrangian coordinates as
\begin{subequations}\label{sys:plf-lg}
    \begin{gather}
        \rho_0 \dt v - \bigl(\mu \dfrac{v_x}{\eta_x}\bigr)_x = 0, \label{eq:plf-lg} \\
        \mu \dx v\vert_{x=-1,1} = 0. \label{bc:sf-lg}
    \end{gather}
\end{subequations}

\bigskip 

We will show that the flow with constant viscosity coefficient $ \mu > 0 $ will become stationary with finite radius. Moreover, we give the characteristic of the asymptotic behavior, including
\begin{itemize}
    \item the terminal domain radius in term of the initial data;
    \item the decay estimates of the velocity; and
    \item global regularity of the solution. 
\end{itemize}

\subsection{Global point-wise estimates}\label{subsec:pointwise-pressureless}

{\noindent\bf Basic energy: }
Taking the $ L^2 $-inner product of \eqref{eq:plf-lg} with $ v $ and applying integration by parts yield
\begin{equation}\label{cnt-vis:001}
    \dfrac{d}{dt} \dfrac{1}{2} \int \rho_0 \vert v \vert^2 \idx + \int \mu \dfrac{\vert v_x \vert^2}{\eta_x} \idx = 0.
\end{equation}
Therefore integrating \eqref{cnt-vis:001} in time implies
\begin{equation}\label{cnt-vis:002}
    \sup_{t} \dfrac{1}{2} \int \rho_0 \vert v \vert^2 \idx + \int_0^\infty \int \mu \dfrac{\vert v_x \vert^2}{\eta_x} \idx \,dt = \dfrac{1}{2} \int \rho_0 \vert v_0 \vert^2 \idx.
\end{equation}

{\noindent\bf Point-wise estimate: }
On the other hand, integrating \eqref{eq:plf-lg} from $ x = -1 $ to $ x \in (-1,1) $ leads to
\begin{equation}\label{cnt-vis:003}
    \dt \int_{-1}^x \rho_0 v \,dx' = \mu \dfrac{v_x}{\eta_x} = \mu \dt \log \eta_x.
\end{equation}
Integrating \eqref{cnt-vis:003} in time yields
\begin{equation}\label{cnt-vis:004}
    \mu \log \eta_x(t) = \int_{-1}^x \rho_0 v(t) \,dx' + \mu \underbrace{\log \eta_{0,x}}_{=0} - \int_{-1}^x \rho_0 v_0 \,dx'.
\end{equation}
Therefore, after applying H\"older's inequality and the energy identity \eqref{cnt-vis:002}, one has
\begin{equation}\label{cnt-vis:005}
    \sup_{x,t} \vert \mu \log \eta_x(x,t) \vert \lesssim \sup_t \int \rho_0 \vert v \vert^2 \idx \leq \int \rho_0 \vert v_0 \vert^2 \idx,
\end{equation}
which implies that, there is a constant $ c_1 \in (0,\infty) $, such that
\begin{equation}\label{cnt-vis:006}
    \dfrac{1}{c_1} \leq \eta_x(x,t) \leq c_1, \qquad (x,t) \in (-1,1)\times [0,\infty).
\end{equation}

{\noindent\bf Decay estimate: }
Thanks to \eqref{cnt-vis:006}, one can conclude from \eqref{cnt-vis:001} that
\begin{equation}\label{cnt-vis:007}
    \dfrac{d}{dt} \int \rho_0 \vert v \vert^2 \idx + \dfrac{2}{c_1} \int \mu \vert v_x \vert^2 \idx \leq 0.
\end{equation}
Meanwhile one has the following Poincar\'e type inequality:
\begin{equation}\label{ineq:poincare-type}
    \int \rho_0 \vert v \vert^2 \idx \lesssim \int \vert v_x \vert^2 \idx. 
\end{equation}
Therefore, \eqref{cnt-vis:007} implies that there exists $ c_2 \in (0,\infty) $ such that
\begin{equation}\label{cnt-vis:008}
    \dfrac{d}{dt} \int \rho_0 \vert v \vert^2 \idx + c_2 \int \rho_0 \vert v \vert^2 \idx \leq 0.
\end{equation}
Therefore, one can conclude that
\begin{equation}\label{cnt-vis:009}
    \int \rho_0 \vert v(t) \vert^2 \idx \leq e^{-c_2 t} \int \rho_0 \vert v_0 \vert^2 \idx. 
\end{equation}
Multiplying \eqref{cnt-vis:007} with $ e^{c_2t/2} $ and integrating the resultant in $ t $ yield
\begin{equation}\label{cnt-vis:011}
    \int_0^\infty e^{c_2 t/2} \norm{v_x(t)}{L^2}^2 \,dt \leq \initial.
\end{equation}

\smallskip

Recalling \eqref{cnt-vis:004}, one can conclude that from \eqref{cnt-vis:009}
\begin{equation}\label{cnt-vis:010}
    \begin{gathered}
        \vert \mu \log \eta_x (x,t) + \int_{-1}^x \rho_0 v_0 \,dx'\vert \leq (\int \rho_0 \idx)^{1/2} (\int \rho_0 \vert v(t) \vert^2 \idx)^{1/2} \\
        \leq e^{-c_2 t/2} (\int \rho_0 \idx)^{1/2} (\int \rho_0 \vert v_0 \vert^2 \idx)^{1/2}.
    \end{gathered}
\end{equation}

\smallskip

Notably, \eqref{cnt-vis:010} states that in the original Euclidean coordinates, the size of the domain (given by $ \int_{-1}^{1}\eta_x \idx $) will be finite all the time, and \eqref{cnt-vis:010} has described the final state of the domain as $ t \rightarrow \infty $. 

\subsection{Global regularity estimates}\label{subsec:regularity-pressureless}

{\noindent\bf $ H^1 $ estimate of $ v $:}

Thanks to \eqref{cnt-vis:003} and \eqref{cnt-vis:006},
one has that  
\begin{equation}\label{cnt-vis:101}
    \norm{v_x(t)}{L^\infty}^2 + \norm{v_x(t)}{L^2}^2 \lesssim \int \rho_0 \vert v_t (t)\vert^2 \,dx 
\end{equation}
for any $ t \in (0,\infty) $.

On the other hand, applying $ \partial_t $ to \eqref{eq:plf-lg} yields 
\begin{equation}\label{cnt-vis:102}
    \rho_0 \partial_t v_t - \bigl( \mu \dfrac{v_{xt}}{\eta_x} \bigr)_x = - \bigl( \mu \dfrac{v_x^2}{\eta_x^2} \bigr)_x.
\end{equation}
Taking the $ L^2 $-inner product of \eqref{cnt-vis:102} with $ v_t $ and applying integration by parts, similar to \eqref{cnt-vis:001}, lead to
\begin{equation}\label{cnt-vis:103}
\begin{aligned}
    & \dfrac{d}{dt}\dfrac{1}{2}\int \rho_0 \vert v_t \vert^2 \,dx + \int \mu \dfrac{\vert v_{xt} \vert^2}{\eta_x} \,dx =  \int \mu \dfrac{v_x^2 v_{xt}}{\eta_x^2} \,dx\\
    & \leq \dfrac{1}{2} \int \mu \dfrac{\vert v_{xt} \vert^2}{\eta_x} \,dx + \dfrac{1}{2} \norm{\frac{v_x}{\eta_x}}{L^\infty}^2 \int \mu \dfrac{\vert v_x \vert^2}{\eta_x} \,dx.
\end{aligned} 
\end{equation}
Consequently, combining \eqref{cnt-vis:006}, \eqref{cnt-vis:101}, and \eqref{cnt-vis:103} implies
\begin{equation}\label{cnt-vis:104}
    \begin{gathered}
    \dfrac{d}{dt} \int \rho_0 \vert v_t \vert^2 \,dx + c_2 \int \rho_0 \vert v_t \vert^2 \,dx 
    \leq \dfrac{d}{dt} \int \rho_0 \vert v_t \vert^2 \,dx + \int \mu \dfrac{\vert v_{xt} \vert^2}{\eta_x} \,dx \\ 
    \lesssim \int \rho_0 \vert v_t \vert^2 \,dx \times \norm{v_x(t)}{L^2}^2,  
    \end{gathered}
\end{equation}
where we have applied the Poincar\'e type inequality \eqref{ineq:poincare-type} with $ v $ replaced by $ v_t $. Therefore, thanks to \eqref{cnt-vis:011}, one can get, after applying Gr\"onwall's inequality in \eqref{cnt-vis:104}, that
\begin{equation}\label{cnt-vis:105}
    \int \rho_0 \vert v_t(t) \vert^2 \,dx \lesssim e^{-c_2t} \times \int \rho_0 \vert v_{\mrm{in},t} \vert^2 \,dx \times e^{C\int_0^t \norm{v_x(s)}{L^2}^2\,ds} \lesssim e^{-c_2 t} \times \initial.
\end{equation}
Therefore, \eqref{cnt-vis:101} and \eqref{cnt-vis:105} implies
\begin{equation}\label{cnt-vis:106}
    \norm{v_x(t)}{L^\infty}^2 + \norm{v_x(t)}{L^2}^2 \lesssim e^{-c_2t}\times \initial.
\end{equation}
Meanwhile, similar to \eqref{cnt-vis:011}, one has that
\begin{equation}\label{cnt-vis:107}
    \int_0^\infty e^{c_2t/2} \norm{v_{xt}(t)}{L^2}^2 \,dt \lesssim \initial.
\end{equation}

{\par\noindent\bf $ H^2 $ estimate of $ v $:} After rewriting \eqref{eq:plf-lg} as 
\begin{equation}\label{cnt-vis:108}
    \mu v_{xx} = \mu \dfrac{v_x}{\eta_x} \eta_{xx} + \rho_0 v_t \eta_x,
\end{equation}
taking the $ L^2 $-inner product of \eqref{cnt-vis:108} with $ \eta_{xx} $ yields
\begin{equation}\label{cnt-vis:109}
    \dfrac{d}{dt} \dfrac{\mu}{2} \norm{\eta_{xx}}{L^2}^2 \leq \mu \norm{\frac{v_x}{\eta_x}}{L^\infty} \norm{\eta_{xx}}{L^2}^2 + \norm{\eta_x}{L^\infty} \bigl( \int \rho_0^2 \vert v_t \vert^2 \,dx \bigr)^{1/2} \norm{\eta_{xx}}{L^2}.
\end{equation}
Therefore, thanks to \eqref{cnt-vis:006}, \eqref{cnt-vis:105}, and \eqref{cnt-vis:106}, one can obtain from \eqref{cnt-vis:109} that
\begin{equation}\label{cnt-vis:110}
    \norm{\eta_{xx}(t)}{L^2} \leq \initial,
\end{equation}
for all $ t \in (0,\infty) $. Using \eqref{cnt-vis:108}, \eqref{cnt-vis:110} together with \eqref{cnt-vis:006}, \eqref{cnt-vis:105}, and \eqref{cnt-vis:106} yields
\begin{equation}\label{cnt-vis:111}
    \norm{v_{xx}(t)}{L^2} \lesssim e^{-c_2t/2} \times \initial,
\end{equation}
for all $ t \in (0,\infty) $. 

\section{Almost pressureless heat conductive flow}\label{sec:almost-pressureless}

Considering system \eqref{sys:NSF-flow}, we define the corresponding Lagrangian coordinates in the same way as in \eqref{def:lagrangian}. Denote by, in addition to \eqref{def:lg-variables}, 
\begin{equation}\label{def:lg-variables-2}
    \Theta(x,t) := \theta(\eta(x,t),t). 
\end{equation}

Then it is straightforward to rewrite system \eqref{sys:NSF-flow} in the lagrangian coordinates as follows:
\begin{subequations}\label{sys:nsf-lg}
    \begin{gather}
        \rho_0 \dt v + \dfrac{1}{M^2} \bigl( \dfrac{\rho_0\Theta}{\eta_x} \bigr)_x - \bigl( \mu \dfrac{v_x}{\eta_x} \bigr)_x = 0, \label{eq:nsf-v-lg}\\
        \rho_0 \dt\Theta + \dfrac{1}{M^2} \rho_0 \Theta \dfrac{v_x}{\eta_x} - \bigl( \kappa \dfrac{\Theta_x}{\eta_x} \bigr)_x = \mu \dfrac{\vert v_x\vert^2}{\eta_x}, \label{eq:nsf-theta-lg} \\
        \mu \partial_x v\vert_{x=-1,1} = \Theta_{x=-1,1} = 0. \label{bc:nsf-lg}
    \end{gather}
\end{subequations}

\subsection{Point-wise estimates}\label{subsec:pointwise-nsf}

{\noindent\bf Basic energy: } Taking the $ L^2 $-inner product of \eqref{eq:nsf-v-lg} with $ v $, integrating \eqref{eq:nsf-theta-lg} in $ I $, and summing up the resultants lead to
\begin{equation}\label{nsf:001}
    \dfrac{d}{dt} \biggl\lbrace \dfrac{1}{2} \int \rho_0 \vert v \vert^2 \idx + \int \rho_0 \Theta \idx \biggr\rbrace = \kappa \dfrac{\Theta_x}{\eta_x} \Big\vert_{x=-1}^{x=1} \leq 0,
\end{equation}
thanks to the fact $ \Theta \geq 0 $ in $ I $. 
Therefore integrating \eqref{nsf:001} in time yields 
\begin{equation}\label{nsf:002}
    \sup_t \biggl\lbrace \dfrac{1}{2} \int \rho_0 \vert v \vert^2 \idx + \int \rho_0 \Theta \idx \biggr\rbrace \leq \biggl\lbrace \dfrac{1}{2} \int \rho_0 \vert v_0 \vert^2 \idx + \int \rho_0 \Theta_0 \idx \biggr\rbrace.
\end{equation}

{\noindent\bf Point-wise estimate: }

Integrating \eqref{eq:nsf-v-lg} from $ x=-1 $ to $ x \in (-1,1) $ leads to
\begin{equation}\label{nsf:003}
    \dt \int_{-1}^x \rho_0 v \,dx' + \dfrac{1}{M^2}\dfrac{\rho_0 \Theta }{\eta_x} = \mu \dfrac{\dt \eta_x}{\eta_x}.
\end{equation}
\eqref{nsf:003} is an ODE of $ \eta_x $ for any fixed $ x \in (-1,1) $. One can solve for $ \eta_x $ as
\begin{equation}\label{nsf:004}
    \begin{gathered}
    \eta_x(t) = \exp \bigl( \frac{1}{\mu}\int_{-1}^x \rho_0 v(t) \idx -  \frac{1}{\mu} \int_{-1}^x \rho_0 v_0 \idx \bigr)\\
    + \frac{1}{\mu M^2}\int_0^t \rho_0\Theta(s) e^{\frac{1}{\mu} \int_{-1}^x \rho_0v(t) \,dx'- \frac{1}{\mu} \int_{-1}^x \rho_0 v(s) \,dx'} \,ds. 
\end{gathered}
\end{equation}
Thanks to the fact $ \Theta \geq 0 $ and $ \int \rho_0 v \,dx \leq \bigl(\int \rho_0 \idx\bigr)^{1/2} \bigl( \int\rho_0 \vert v \vert^2 \idx \bigr)^{1/2} $, \eqref{nsf:002} and \eqref{nsf:004} imply that 
\begin{equation}\label{nsf:005}
    0 < \frac{1}{c_3} \leq \eta_x \leq c_3 + \frac{c_3}{\mu M^2} \int_0^t \rho_0 \Theta(s) \,ds,
\end{equation}
for some $ c_3 \in (0,\infty)$. Notice that, unlike the pressureless flow \eqref{cnt-vis:006}, there is no {\it a priori} bound of the evolving boundary from the basic energy estimate \eqref{nsf:005}.

To continue, we assume that, for $ M $ large enough, 
\begin{equation}\label{asmp:aprior}
    \frac{1}{\mu M^2} \int_0^t \rho_0 \Theta(s) \,ds \leq 1, 
\end{equation}
hence 
\begin{equation}\label{asmp:aprior-2}
    0 < \frac{1}{2 c_3} \leq \frac{1}{\eta_x} \leq 2 c_3 
\end{equation}

\subsection{Estimates under {\it a prior} assumption \eqref{asmp:aprior}}
\label{subsec:estimates-nsf}


{\noindent\bf Estimates on the velocity $ v $.} Assuming \eqref{asmp:aprior-2}, one can derive from \eqref{eq:nsf-v-lg} that, after taking the $ L^2 $-inner product with $ 2 v $ and integration by parts,
\begin{equation}\label{nsf:006}
    \begin{gathered}
        \dfrac{d}{dt} \int \rho_0 \vert v \vert^2 \idx + \frac{1}{c_3} \int  \mu \vert v_x \vert^2 \idx \leq \frac{1}{M^2} \int \dfrac{\rho_0 \Theta v_x}{\eta_x} \idx \\
        \leq \frac{\mu}{2c_3} \int \vert v_x \vert^2 \idx + \frac{c_3}{2\mu M^4} \int \dfrac{\rho_0^2 \Theta^2}{\eta_x^2} \idx.
\end{gathered}
\end{equation}
Therefore one has that, after applying \eqref{ineq:poincare-type}, 
\begin{equation}\label{nsf:007}
    \dfrac{d}{dt} \int \rho_0 \vert v \vert^2 \idx + c_4 \int \rho_0 \vert v \vert^2 \idx + c_4 \int \vert v_x \vert^2 \idx \leq \frac{c_5}{M^4} \int \rho_0 \Theta^2 \idx,
\end{equation}
for some constants $ c_4, c_5 $ independent of $ M $. Therefore, after applying Gr\"onwall's inequality, \eqref{nsf:007} implies 
\begin{equation}
    \label{nsf:0071}
    \begin{gathered}
        \int \rho_0 \vert v(t) \vert^2 \idx + c_4 e^{-c_4t} \int_0^t e^{c_4 s} \norm{v_x(s)}{L^2}^2 \idx\,ds \leq e^{-c_4t} \int \rho_0 \vert v_0 \vert^2 \idx \\
        + \frac{c_5}{M^4} \int_0^t e^{c_4(s-t)} \int \rho_0 \Theta^2(s) \idx\,ds.
    \end{gathered}
\end{equation}       

\smallskip 

Meanwhile, \eqref{nsf:003} implies
\begin{gather}\label{nsf:008}
    \norm{\frac{v_x}{\eta_x}}{L^\infty}^2 \lesssim  \int \rho_0 \vert v_t \vert^2 \idx + \dfrac{1}{M^4} \norm{\Theta}{\infty}^2 \lesssim \int \rho_0 \vert v_t \vert^2 \idx + \dfrac{1}{M^4} \norm{\Theta_x}{L^2}^2,   \\
    \intertext{and}
    \label{nsf:0081}
    \norm{\frac{v_x}{\eta_x}}{L^2}^2 \lesssim \int \rho_0 \vert v_t \vert^2 \idx + \dfrac{1}{M^4} \int \rho_0 \vert\Theta \vert^2 \idx,
\end{gather}
thanks to \eqref{asmp:aprior-2}, where the Poincar\'e inequality
\begin{equation}\label{nsf:0082}
    \norm{\Theta}{L^\infty} + \norm{\Theta}{L^2} \lesssim \norm{\Theta_x}{L^2}. 
\end{equation}
has applied. 

Applying $ \partial_t $ to \eqref{eq:nsf-v-lg} yields
\begin{equation}\label{nsf:009}
        \rho_0 \dt v_t - \bigl( \mu \dfrac{v_{xt}}{\eta_x} \bigr)_x =  - \bigl( \mu \dfrac{v_x^2}{\eta_x^2} \bigr)_x - \frac{1}{M^2} \bigl( \dfrac{\rho_0 \Theta}{\eta_x} \bigr)_{xt}.
\end{equation}
Taking the $ L^2 $-inner product of \eqref{nsf:009}  with $ v_t $ and applying integration by parts lead to
\begin{equation}\label{nsf:010}
\begin{gathered}
    \dfrac{d}{dt} \frac{1}{2} \int \rho_0 \vert v_t \vert^2 \idx + \int \mu \dfrac{\vert v_{xt} \vert^2}{\eta_x} \idx = \int \mu \dfrac{v_x^2 v_{xt}}{\eta_x^2} + \dfrac{1}{M^2} \int  \bigl( \dfrac{\rho_0 \Theta}{\eta_x} \bigr)_t v_{xt} \idx \\
    \leq \frac{\mu}{2} \int \dfrac{\vert v_{xt}\vert^2}{\eta_x} \idx + \frac{1}{2} \norm{\dfrac{v_x}{\eta_x}}{L^\infty}^2 \int \mu \dfrac{v_x^2}{\eta_x} \idx + \frac{1}{2\mu M^4} \int \eta_x \vert \bigl(\dfrac{\rho_0 \Theta}{\eta_x} \bigr)_t \vert^2 \idx.
\end{gathered}
\end{equation}
Combining \eqref{asmp:aprior-2}, \eqref{nsf:008}, \eqref{nsf:0082}, \eqref{nsf:010} and the Poincar\'e type inequality \eqref{ineq:poincare-type} for $ v_t $ implies
\begin{equation}
    \label{nsf:011}
    \begin{gathered}
        \dfrac{d}{dt} \int \rho_0 \vert v_t \vert^2 \idx + c_4 \int \rho_0 \vert v_t \vert^2 \idx + c_4 \norm{v_{xt}}{L^2}^2 \\
        \leq c_6  \norm{v_x}{L^2}^2 \bigl( \int\rho_0 \vert v_t \vert^2 \idx + \frac{1}{M^4}\norm{\Theta_x}{L^2}^2 \bigr) \\
        + \frac{c_6}{M^4} \bigl( \int \rho_0 \Theta_t^2 \idx + \norm{\Theta_x}{L^2}^2 \norm{v_x}{L^2}^2 \bigr),
    \end{gathered} 
\end{equation}
for some positive constant $ c_6 $. Here, one can carefully check that $ c_4 $ is the same as in \eqref{nsf:007}. Then, applying Gr\"onwall's inequality in \eqref{nsf:011} yields
\begin{equation}
    \label{nsf:012}
\begin{gathered}
    e^{c_4 t - c_6 \int_0^t\norm{v_x(s)}{L^2}^2\,ds}\int \rho_0 \vert v_t(t) \vert^2 \idx + c_4 \int_0^t e^{c_4 s - c_6 \int_0^s\norm{v_x(\tau)}{L^2}^2\,d\tau} \norm{v_{xt}(s)}{L^2}^2 \,ds \\
    \leq \int \rho_0\vert v_1 \vert^2 \idx \\
    + \frac{2c_6}{M^4} \int_0^t e^{c_4 s - c_6\int_0^s \norm{v_x(\tau)}{L^2}^2 \,d\tau} (\norm{v_x(s)}{L^2}^2 \norm{\Theta_x(s)}{L^2}^2 + \int \rho_0 \vert \Theta_t(s) \vert^2 \idx)\,ds.
\end{gathered}
\end{equation}

{\bigskip \noindent\bf Estimates on the temperature $ \Theta $.} Taking the $ L^2 $-inner product of \eqref{eq:nsf-theta-lg} with $ \Theta $ and applying integration by parts in the resultant lead to 
\begin{equation}
    \label{nsf:013}
    \begin{gathered}
        \dfrac{d}{dt} \frac{1}{2} \int \rho_0 \vert \Theta \vert^2 \idx  + \kappa \int \dfrac{\vert \Theta_x \vert^2}{\eta_x} \idx = \mu \int \dfrac{\vert v_x \vert^2 \Theta}{\eta_x} \idx \\
        - \frac{1}{M^2} \int \dfrac{\rho_0 \vert\Theta\vert^2 v_x}{\eta_x} \idx \lesssim \norm{\Theta_x}{L^2} \norm{v_x}{L^2}^2 + \frac{1}{M^2} \norm{\Theta_x}{L^2} \norm{v_x}{L^2} \norm{\Theta_x}{L^2},
    \end{gathered}
\end{equation}
where the Poincar\'e type inequality \eqref{nsf:0082} has applied. Then applying the Cauchy-Schwarz inequality in \eqref{nsf:013} implies that, there exists $ c_7, c_8 >0 $ that
\begin{equation}
    \label{nsf:014}
    \dfrac{d}{dt} \int \rho_0 \vert \Theta \vert^2 \idx + c_7 \int \rho_0 \vert \Theta \vert^2 \idx + c_7 \norm{\Theta_x}{L^2}^2 \leq c_8 \norm{v_x}{L^2}^4 + \frac{c_8}{M^4} \norm{v_x}{L^2}^2 \norm{\Theta_x}{L^2}^2, 
\end{equation} 
where the Poincar\'e type inequality \eqref{nsf:0082} has applied.
Then applying Gr\"onwall's inequality in \eqref{nsf:014} implies
\begin{equation}
    \label{nsf:015}
\begin{gathered}
    e^{c_7 t} \int \rho_0 \vert\Theta(t)\vert^2 \idx  + c_7 \int_0^t e^{c_7s} \norm{\Theta_x(s)}{L^2}^2 \,ds \leq \int \rho_0 \vert \Theta_0 \vert^2 \idx \\
    + c_8 \int_0^t e^{c_7s} \norm{v_x(s)}{L^2}^4 \,ds + \frac{c_8}{M^4} \int_0^t e^{c_7s} \norm{v_x(s)}{L^2}^2 \norm{\Theta_x(s)}{L^2}^2 \,ds.
\end{gathered}
\end{equation}

Meanwhile, applying $ \partial_t $ to \eqref{eq:nsf-theta-lg} yields
\begin{equation}
    \label{nsf:016}
    \begin{gathered}
        \rho_0 \partial_t \Theta_t - \bigl( \kappa \dfrac{\Theta_{xt}}{\eta_x} \bigr)_x = - \bigl( \kappa \dfrac{\Theta_x v_x}{\eta_x^2} \bigr)_x + \mu \dfrac{2 v_x v_{xt}}{\eta_x} - \mu \dfrac{v_x^3}{\eta_x^2} \\
        - \frac{1}{M^2} \biggl( \dfrac{\rho_0 \Theta_t v_x}{\eta_x} + \dfrac{\rho_0 \Theta v_{xt}}{\eta_x} - \dfrac{\rho_0 \Theta v_{x}^2}{\eta_x^2} \biggr).  
    \end{gathered}
\end{equation}
After taking the $ L^2 $-inner product of \eqref{nsf:016} with $ \Theta_t $ and applying integration by parts in the resultant, one obtains
\begin{equation}
    \label{nsf:017}
    \begin{gathered}
        \dfrac{d}{dt} \frac{1}{2}\int \rho_0 \vert \Theta_t \vert^2 \idx + \kappa \int \dfrac{\vert \Theta_{xt}\vert^2}{\eta_x} \idx = \kappa \int \dfrac{\Theta_x v_x \Theta_{xt}}{\eta_x^2} \idx \\
        + \mu \int \biggl( \dfrac{2 v_x v_{xt} \Theta_t}{\eta_x} - \dfrac{v_x^3\Theta_t}{\eta_x^2} \biggr) \idx \\
        - \frac{1}{M^2} \int \biggl( \dfrac{\rho_0 \Theta_t^2 v_x}{\eta_x} + \dfrac{\rho_0 \Theta v_{xt}\Theta_t}{\eta_x} - \dfrac{\rho_0 \Theta \Theta_t v_{x}^2}{\eta_x^2} \biggr)\idx\\
        \lesssim \norm{\Theta_{xt}}{L^2} \norm{\Theta_{x}}{L^2} \norm{{v_x}}{L^\infty} + \norm{v_x}{L^2} \norm{v_{xt}}{L^2} \norm{\Theta_{xt}}{L^2} \\
        + \norm{v_x}{L^2}^2 \norm{{v_x}}{L^\infty} \norm{\Theta_{xt}}{L^2} + \frac{1}{M^2} \norm{\Theta_{xt}}{L^2}\norm{v_x}{L^2}  \biggl(\int \rho_0 \vert \Theta_t \vert^2 \idx \biggr)^{1/2} \\
        + \dfrac{1}{M^2} \biggl( \norm{v_{xt}}{L^2} \norm{\Theta_{xt}}{L^2} + \norm{\Theta_{xt}}{L^2} \norm{v_x}{L^\infty} \norm{v_x}{L^2}  \biggr) \biggl( \int \rho_0 \vert \Theta \vert^2 \idx  \biggr)^{1/2}
    \end{gathered}
\end{equation}
where we have applied the Poincar\'e inequality \eqref{nsf:0082} for $ \Theta_t $.
Then similarly to \eqref{nsf:014}, after applying Cauchy-Schwarz inequality in \eqref{nsf:017}, one has that, for some $ c_9 >0 $ that
\begin{equation}
    \label{nsf:018}
    \begin{gathered}
        \dfrac{d}{dt} \int \rho_0 \vert \Theta_t \vert^2 \idx + c_7 \int \rho_0 \vert \Theta_t \vert^2 \idx + c_7 \norm{\Theta_{xt}}{L^2}^2 \\
        \leq c_9 \bigl( \norm{\Theta_x}{L^2}^2 \norm{v_x}{L^\infty}^2 + \norm{v_{xt}}{L^2}^2 \norm{v_x}{L^2}^2 + \norm{v_x}{L^2}^4 \norm{v_x}{L^\infty}^2 \bigr) \\
        + \frac{c_9}{M^4} \biggl( \norm{v_x}{L^2}^2 \int\rho_0 \vert \Theta_t \vert^2 \idx + \bigl( \norm{v_{xt}}{L^2}^2 + \norm{v_x}{L^2}^2 \norm{v_x}{L^\infty}^2 \bigr) \int \rho_0 \vert \Theta \vert^2 \idx  \biggr).
    \end{gathered}
\end{equation}
Applying Gr\"onwall's inequality in \eqref{nsf:018} yields
\begin{equation}
    \label{nsf:019}
    \begin{gathered}
        e^{c_7 t} \int \rho_0 \vert \Theta_t (t) \vert^2 \idx + c_7 \int_0^t e^{c_7 s} \norm{\Theta_{xt}(s)}{L^2}^2  \,ds \leq \int \rho_0 \vert \Theta_1 \vert^2 \idx  \\
        + c_9 \int_0^t e^{c_7 s} \bigl( \norm{\Theta_x}{L^2}^2 \norm{v_x}{L^\infty}^2 + \norm{v_{xt}}{L^2}^2 \norm{v_x}{L^2}^2 + \norm{v_x}{L^2}^4 \norm{v_x}{L^\infty}^2 \bigr) \,ds \\
        + \frac{c_9}{M^4} \int_0^t e^{c_7 s} \biggl( \substack{\norm{v_x}{L^2}^2 \int\rho_0 \vert \Theta_t \vert^2 \idx + \bigl( \norm{v_{xt}}{L^2}^2 + \norm{v_x}{L^2}^2 \norm{v_x}{L^\infty}^2 \bigr) \\
        \times  \int \rho_0 \vert \Theta \vert^2 \idx}  \biggr) \,ds.
    \end{gathered}
\end{equation}

\subsection{Collecting the {\it a prior} estimates}
\label{subsec:collecting-nsf}

From \eqref{nsf:007}, \eqref{nsf:011}, \eqref{nsf:014}, and \eqref{nsf:018}, by choosing the smallest one among $ c_4, c_7 $ and largest among $ c_5, c_6, c_8, c_9 $, we can without loss of generality assume that
\begin{equation}\label{nsf:020}
    \mathfrak{c}_1 = c_4 = c_7 \in (0,\infty), \qquad \mathfrak{c}_2 = c_5 = c_6 = c_8 = c_9 \in (0,\infty). 
\end{equation}
We also write $ \mathfrak c_3 \in (0,\infty) $ in the following, to represent a large but finite constant which may be different from line to line. 
Furthermore, denote by
\begin{align}
    & \label{def:nsf-energy}
        \begin{aligned}
            \text{Energy:}  \quad \mathcal E(t) := & e^{\mathfrak c_1 t} \int \rho_0 \vert v (t) \vert^2 \idx + e^{\mathfrak{c}_1} \int \rho_0 \vert v_t(t) \vert^2 \idx \\
            & + e^{\mathfrak c_1 t} \int \rho_0 \vert \Theta (t) \vert^2 \idx + e^{\mathfrak{c}_1} \int \rho_0 \vert \Theta_t(t) \vert^2 \idx\\
            & + \int_0^t \bigl( e^{\mathfrak{c}_1 s} \norm{v_x(s)}{L^2}^2 + e^{\mathfrak c_1 s} \norm{v_{xt}(s)}{L^2}^2 \bigr) \,ds \\
            & + \int_0^t \bigl( e^{\mathfrak{c}_1 s} \norm{\Theta_x(s)}{L^2}^2 + e^{\mathfrak c_1 s} \norm{\Theta_{xt}(s)}{L^2}^2 \bigr) \,ds, 
        \end{aligned}\\
    & \label{def:nsf-initial-energy}
    \begin{aligned} 
        \text{Initial data:} \quad \mathcal E_0 := & \int\rho_0 \vert v_0 \vert^2 \idx + \int \rho_0 \vert v_1 \vert^2 \idx \\ 
        & + \int\rho_0 \vert \Theta_0 \vert^2 \idx + \int \rho_0 \vert \Theta_1 \vert^2 \idx.
    \end{aligned}
\end{align}
Then thanks to \eqref{nsf:008} and \eqref{nsf:0081}, one has
\begin{gather}
    \label{nsf:021} \begin{gathered} e^{\mathfrak c_1 t} \norm{v_x(t)}{L^2}^2 + \int_0^t e^{\mathfrak c_1(s)}(\norm{v_x(s)}{L^2}^2 + \norm{v_x(s)}{L^\infty}^2) \,ds\\
         +\int_0^t e^{\mathfrak c_1 s} \int \rho_0 \vert \Theta_t(s) \vert^2 \idx \,ds  \lesssim  \mathcal E(t), \end{gathered} 
\end{gather} 
where we have applied \eqref{ineq:poincare-type} for $ v_t $ and \eqref{nsf:0082} for $ \Theta_t $. 
Moreover, we use, in the following,  
\begin{equation}
    \label{nsf:0211}
    P(\mathcal E(t))
\end{equation}
to represent a function of $ \mathcal E(t) $ with $ P(0) = 0 $, which will be different from line to line. 

\smallskip
Then from \eqref{nsf:0071}, one has
\begin{equation}
    \label{nsf:023}
    \begin{gathered}
        e^{\mathfrak c_1 t} \int \rho_0 \vert v(t) \vert^2 \idx + \int_0^t e^{\mathfrak c_1 s} \norm{v_x(s)}{L^2}^2 \,ds \lesssim \mathcal E_0 + \frac{1}{M^4} \mathcal E(t).
    \end{gathered}
\end{equation}
With \eqref{nsf:023}, \eqref{nsf:012} implies that,
\begin{equation}
    \label{nsf:024}
    \begin{gathered}
        e^{\mathfrak{c}_1 t} \int \rho_0 \vert v_t(t) \vert^2 \idx + \int_0^t e^{\mathfrak{c}_1 s} \norm{v_{v_{xt}(s)}}{L^2}^2 \,ds \lesssim \mathcal E_0 e^{\mathfrak c_3 \mathcal E_0 + \frac{\mathfrak c_3}{M^4} \mathcal E(t)} \\
        + \frac{e^{\mathfrak c_3 \mathcal E_0 + \frac{\mathfrak c_3}{M^4} \mathcal E(t)}}{M^4} \bigl\lbrack \mathcal E^2(t) + \mathcal E(t) \bigr\rbrack,
    \end{gathered}
\end{equation}
where we have applied \eqref{nsf:021}. In particular, together with \eqref{nsf:008} and \eqref{nsf:0081}, \eqref{nsf:023} and \eqref{nsf:024} imply that
\begin{equation}
    \label{nsf:025}
    \begin{gathered}
        e^{\mathfrak c_1 t} (\norm{v_x(t)}{L^2}^2 + \norm{v_{x}(t)}{L^\infty}^2) + \int_0^t e^{\mathfrak c_1 s} (\norm{v_x(s)}{L^2}^2 + \norm{v_x(s)}{L^\infty}^2 )\, ds \\
        \lesssim \mathcal E_0 + \frac{1}{M^4} P(\mathcal E(t)), 
    \end{gathered}
\end{equation}
where we have applied the following inequality
\begin{equation}
    \label{nsf:026}
    e^{\mathfrak c_1 t} \norm{\Theta_{x}(t)}{L^2}^2 \lesssim \int_0^t e^{\mathfrak c_1 s} (\norm{\Theta_{x}(s)}{L^2}^2 + \norm{\Theta_{xt}(s)}{L^2}^2)\,ds.
\end{equation}

Similarly, \eqref{nsf:015} implies 
\begin{equation}
    \label{nsf:027}
    \begin{gathered}
        e^{\mathfrak c_1 t} \int \rho_0 \vert \Theta (t) \vert^2 \idx + \int_0^t e^{\mathfrak c_1 s} \norm{\Theta_x(s)}{L^2}^2 \,ds \lesssim \mathcal E_0 + \mathcal E_0^2 \\
        + \frac{1}{M^4} P(\mathcal E(t)),
    \end{gathered}
\end{equation}
where we have applied \eqref{nsf:025}.

Finally, applying \eqref{nsf:024}, \eqref{nsf:025}, and \eqref{nsf:027} on the right hand side of \eqref{nsf:019} leads to 
\begin{equation}
    \label{nsf:028}
    \begin{gathered}
        e^{\mathfrak c_1 t} \int \rho_0 \vert \Theta_t(t) \vert^2 \idx + \int_0^t e^{\mathfrak c_1 s} \norm{\Theta_{xt}(s)}{L^2}^2 \,ds \lesssim \mathcal E_0 + \mathcal E_0^3 +  \frac{1}{M^4} P(\mathcal E(t)).
    \end{gathered}
\end{equation}

\bigskip

In conclusion, \eqref{nsf:023}--\eqref{nsf:028} yield that
\begin{equation}
    \label{nsf:029}
    \mathcal E(t) \leq \mathfrak{c}_3 (\mathcal E_0 + \mathcal E_0^3) + \frac{1}{M^4} P(\mathcal E(t)).
\end{equation}
Therefore, for $ M $ large enough, one can conclude from \eqref{nsf:029} that
\begin{equation}
    \label{nsf:030}
    \mathcal E(t) \lesssim \mathcal E_0 + \mathcal E_0^3,
\end{equation}
and together with \eqref{nsf:025} and \eqref{nsf:026}, we have
\begin{equation}
    \label{nsf:031}
    e^{\mathfrak c_1 t} (\norm{v_x(t)}{L^\infty}^2 + \norm{\Theta_{x}(t)}{L^2}^2) + \int_0^t e^{\mathfrak c_1 t} \norm{v_x(s)}{L^\infty}^2 \,ds \lesssim \mathcal E_0 + \mathcal E_0^3.
\end{equation}

\smallskip 

Now we are ready to close the {\it a prior} assumption \eqref{asmp:aprior}. Indeed, thanks to \eqref{nsf:0082} and \eqref{nsf:031}, one can calculate 
\begin{equation}
    \label{nsf:032} 
    \frac{1}{\mu M^2} \int_0^t \rho_0\Theta(s) \,ds \lesssim \frac{1}{M^2} \int_0^t \norm{\Theta_x(s)}{L^2} \,ds \lesssim \frac{\mathcal E_0 + \mathcal E_0^3}{M^2} \int_0^t e^{\mathfrak{c}_1 s / 2} \,ds.
\end{equation}
This verifies \eqref{asmp:aprior} for $ M $ large enough. 

\section{Appendix: Viscous pressureless flow with degenerate viscosity}

Here we would like to demonstrate the effect of degenerate viscosity by constructing self-similar solutions. We consider system \eqref{sys:pressureless-flows} with $ \mu = \mu(\rho) $, i.e., when the viscosity coefficient depends on the density. In particular, $ \mu \rightarrow 0 $ as $ \rho \rightarrow 0 $ on the moving boundary. Notice that, the lagrangian formulation \eqref{sys:plf-lg} still holds, with $ \mu = \mu(\rho_0/\eta_x) $. 

\smallskip

Consider a solution to \eqref{eq:plf-lg} with $ \mu = \rho^\alpha = (\rho_0/\eta_x)^\alpha $, $ \alpha > 0 $ of the form
\begin{equation}\label{dd-vis:001}
    \eta(x,t) = \sigma(t) x. 
\end{equation}
Then \eqref{eq:plf-lg} is reduced to
\begin{equation}\label{dd-vis:002}
    x \rho_0 \sigma''(t) - (\rho_0^\alpha)_x \dfrac{\sigma'(t)}{\sigma^{\alpha+1}(t)} = 0.
\end{equation}

Consider
\begin{equation}\label{dd-vis:003}
    x\rho_0 = - \dfrac{1}{\alpha} (\rho_0^\alpha)_x. 
\end{equation}
That is
\begin{itemize}
    \item when $ \alpha > 1 $, $ (\alpha-1)(1-x^2)/2 = \rho_0^{\alpha-1}(x), x\in (-1,1) $;
    \item when $ \alpha = 1 $, $\rho_0 (x) = \rho_0(0) e^{-x^2/2}, x \in \mathbb R$;
    \item when $ 0 < \alpha < 1 $, $ \rho_0(x) = \dfrac{1}{[(1-\alpha)x^2/2 + \rho_0^{\alpha-1}(0)]^{1/(1-\alpha)}}, x \in \mathbb R$.
\end{itemize}
Then from \eqref{dd-vis:002}, one has
\begin{equation}\label{dd-vis:004}
    \sigma''(t) - (\sigma^{-\alpha})'(t) = 0.
\end{equation}
Integrating \eqref{dd-vis:004} in time yields
\begin{equation}\label{dd-vis:005}
    \sigma'(t) = \sigma^{-\alpha}(t) + \sigma'(0) - \sigma^{-\alpha}(0).
\end{equation}

Now, let 
\begin{equation}\label{dd-vis:006}
    \Gamma := \sigma'(0) - \sigma^{-\alpha}(0).
\end{equation}
\begin{itemize}
    \item If $ \Gamma \geq 0 $, one can conclude from \eqref{dd-vis:005} that $ \lim_{t\rightarrow \infty} \sigma(t) = \infty $ and $ \lim_{t\rightarrow \infty} \sigma'(t) = \Gamma $.
    \item If $ \Gamma < 0 $, one can conclude from \eqref{dd-vis:005} that $ \lim_{t\rightarrow \infty} \sigma(t) = (-\Gamma)^{-1/\alpha} \in (0,\infty) $ and $ \lim_{t\rightarrow \infty} \sigma'(t) = 0 $.
\end{itemize}
We call the case when $ \Gamma \geq 0 $ the large energy self-similar solution, and the case when $ \Gamma < 0 $ the small energy self-similar solution.



\smallskip 

The large energy and the small energy solutions indicate that, the asymptotic behavior of the solutions depends strongly on the structure of the initial data, regardless of the value of $ \alpha > 0 $. Consequently, it is not realistic to investigate the general solution to \eqref{sys:plf-lg} with density depending viscosity without imposing restriction on the initial data. This is very different from the flow with constant viscosity as shown in theorem \ref{thm:pressureless}.

\bibliographystyle{plain}


\end{document}